\documentclass[12pt]{amsart}
\usepackage{amssymb, colordvi,multirow, graphicx, multicol}

\numberwithin{equation}{section}
\newtheorem{theorem}{Theorem}
\numberwithin{theorem}{section}
\newtheorem{proposition}[theorem]{Proposition}

\newtheorem{ex}[theorem]{Example}
\newtheorem{rem}[theorem]{Remark}

\newenvironment{example}{\begin{ex}\rm}{\end{ex}}
\newenvironment{remark}{\begin{rem}\rm}{\end{rem}}
\newcounter{FNC}[page]
\def\fauxfootnote#1{{\addtocounter{FNC}{2}$^\fnsymbol{FNC}$%
     \let\thefootnote\relax\footnotetext{$^\fnsymbol{FNC}$#1}}}

\newcommand\Gal{\mbox{\rm Gal}}

\newcommand{\C}{{\mathbb{C}}}

\newcommand{\calE}{{\mathcal{E}}}

\copyrightinfo{}{}

\title[Galois groups of Schubert problems]{Galois groups of Schubert problems\\
   via homotopy computation}

\author{Anton Leykin}
\address{Institute for Mathematics and its Applications\\
         University of Minnesota\\
         Minneapolis\\
         MN \ 55455\\
         USA}
\email{leykin@ima.umn.edu}
\urladdr{http://www.ima.umn.edu/~leykin/}

\author{Frank Sottile}
\address{Department of Mathematics\\
         Texas A\&M University\\
         College Station\\
         Texas \ 77843\\
         USA}
\email{sottile@math.tamu.edu}
\urladdr{http://www.math.tamu.edu/\~{}sottile/}

\thanks{Leykin and Sottile supported by the Institute for Mathematics and its Applications and
  Sottile by NSF grants CAREER DMS-0538734 and  DMS-0701050}

\begin{document}

\begin{abstract}
 Numerical homotopy continuation of solutions to polynomial equations is the foundation
 for numerical algebraic geometry, whose development has been driven by
 applications of mathematics.
 We use numerical homotopy continuation to investigate the problem
 in pure mathematics of determining Galois groups in the Schubert calculus.
 For example, we show by direct computation that the Galois group of the Schubert problem
 of 3-planes in $\C^8$ meeting 15 fixed 5-planes non-trivially is the full symmetric group
 $S_{6006}$.
\end{abstract}

\subjclass[2000]{14N15, 65H20}
%
%
\keywords{Polynomial homotopy continuation, Schubert problem, Galois group}
\maketitle
%
\section*{Introduction}
%

Numerical homotopy continuation~\cite{SW05} gives a method to find all solutions to
a system of polynomials with finitely many solutions.
Current parallel implementations~\cite{LVY06}
can solve systems with over~$40$ million solutions~\cite{TY}.
The emerging field of \Blue{{\it numerical algebraic geometry}}~\cite{SVW9,SW05}
uses numerical homotopy continuation as a foundation for algorithms to study
algebraic varieties.
While numerical algebraic geometry was developed for applications of mathematics,
we apply it in pure mathematics, computing Galois groups of
enumerative-geometric problems from the Schubert calculus, called
\Blue{{\it Schubert problems}}.

Along with~\cite{BPSW}, this is one of the first applications of
numerical algebraic geometry
to a problem in pure mathematics.

Jordan introduced Galois groups of enumerative problems
in 1870~\cite{J1870} and
Harris laid their modern foundations in 1979~\cite{Ha79},
showing that the algebraic Galois group is equal to
a geometric monodromy group.
Byrnes~\cite[Section 5]{By89} used Harris's theory to prove that the general problem of
placing poles with static output feedback in linear systems theory was not solvable by
radicals.
He used numerical homotopy continuation to show that a particular
Galois group arising in pole placement was the full symmetric group,~$S_5$.
Underlying this calculation was a Schubert problem.
Vakil~\cite{Va06b} applied his geometric Littlewood-Richardson rule~\cite{Va06a}
to study Galois groups of Schubert problems and showed
that many Schubert problems have Galois group containing the alternating group.

A Schubert problem is \Blue{{\it simple}} if it involves no more than two Schubert
conditions of codimension more than 1.
Simple Schubert problems are natural to study~\cite{So00b,So03} and among all
Schubert problems on a given Grassmannian, they have the largest intersection numbers,
so they are the most challenging for direct computation.
They may also be formulated as complete intersections, which is a restriction imposed by
our software.
\medskip

\noindent{\bf Numerical Theorem.}
{\it
  The Galois group of the Schubert problem of\/ $3$-planes in $\C^8$
  meeting $15$ fixed $5$-planes non-trivially is the full symmetric group
  $S_{6006}$. \vspace{10pt}
}

This is a numerical theorem, as our software does not certify its output.
We have computed Galois groups of scores of other simple Schubert
problems, including one on the Grassmannian of $4$-planes in $\C^8$ having 8580
solutions and ones on  Grassmannians of $3$-planes in $\C^8$ and in $\C^9$ having 10329
and 17589 solutions, respectively. 
In every case, we find that the Galois group is the full symmetric group.
Table~\ref{table:timings} in Section~\ref{Sec:software} records some of these
calculations.
Based on this evidence, we conjecture that every simple Schubert problem on a
Grassmannian has Galois group equal to the full symmetric group.

Not all Schubert problems have Galois group equal to the full symmetric group.
Using an idea of Derksen, Vakil~\cite{Va06b} gives some Schubert problems on Grassmannians
whose Galois group is not the full symmetric group and Ruffo, et.~al.~\cite{RSSS} give one
on a particularly small flag manifold.
None of these examples can be studied with our software, which requires the
Schubert problem to have a formulation as a complete intersection.

Our software has two implementations in Maple
which use homotopy continuation to compute elements in the Galois groups
and either Maple or GAP~\cite{GAP4} to determine if these elements generate the full
symmetric group.
For the continuation, both implementations use PHCpack~\cite{V99}
through its Maple interface PHCmaple~\cite{PHCmaple-ICMS-06} and the second may also call
Bertini~\cite{Bertini}.
The advantages of Bertini are that it can use arbitrary precision and it gives an independent
verification of our results.

Numerical techniques give insight into some mathematical properties that are far beyond
the reach of other methods.
For example, Billey and Vakil~\cite{BV05} studied Galois groups of Schubert problems
using symbolic methods.
The largest problem that they treated (showing its Galois group is the full
symmetric group) had 9 solutions on the Grassmannian of 2-planes in $\C^6$, and
they stated that the Schubert problem on this Grassmannian having 14 solutions
was computationally infeasible.

The largest simple Schubert problem which we have solved
symbolically has 91 solutions~\cite[\S 5.3]{RSSS}.
In contrast, numerical methods allow us to solve Schubert problems with as many as 17589
solutions. 
These examples actually underestimate the gap between the computational possibilities of
symbolic and numeric methods, because they were performed on serial machines.

Current and (likely) future increases in computer
power will come from multiple core and distributed computing.
This is a break with the past, when improvements in computational power came from 
increasing the clock speed of single-processor units.
Symbolic algorithms have limited potential in this regime, as
Gr\"obner basis computation appears to be intrinsically serial and thus can not
be efficiently parallelized.
In contrast, numerical homotopy continuation is easily parallelized,
since its atomic tasks are independent.
Thus methods based on numerical continuation will reap the benefits of future
parallel architectures.
In addition, numerical algorithms typically require less memory than symbolic algorithms.
In particular, the sizes of final and intermediate expressions in Gr\"obner basis
computation not only may be large, but also are unpredictable.
For these reasons, we feel that the future of computing in algebraic geometry lies in
numerical algorithms.

This paper is structured as follows.
In Section~\ref{Sec:Schubert}, we describe the basic geometry of Schubert problems
and Harris's theory of Galois groups.
In Section~\ref{Sec:Homotopy}, we explain the use of homotopy continuation for simple
Schubert problems.
We present our software and algorithms and discuss our results in Section 3,
which include the computation described in the Numerical Theorem.
We describe future work in Section 4.

%
\section{Galois group computation of Schubert problems}\label{Sec:Schubert}
%

\begin{figure}[htb]
\[
  \begin{picture}(314,192)
   \put(3,0){\includegraphics[height=6.6cm]{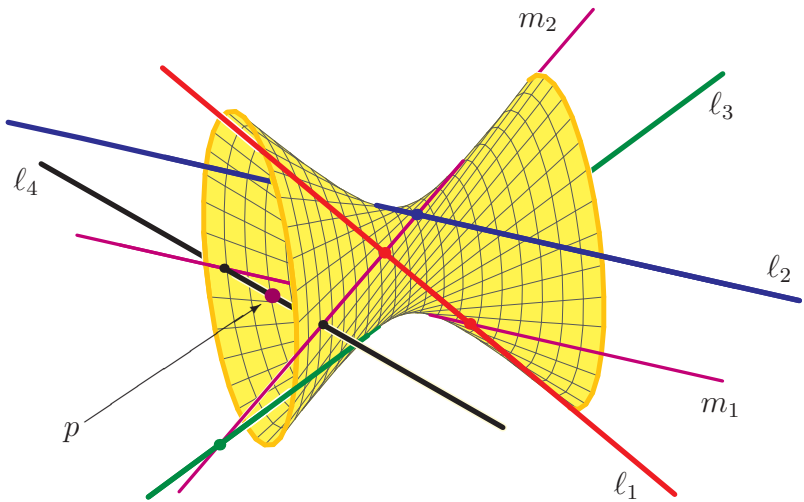}}
   \put(230,  3){$\ell_1$}
   \put(288, 85){$\ell_2$}
   \put(266,147){$\ell_3$}
   \put(  3,118){$\ell_4$}
   \put(263, 34){$m_1$}
   \put(194,180){$m_2$}

   \put( 22, 25){$p$}\put(30,30.5){\vector(3,2){66.6}}
  \end{picture}
\]
\caption{The two lines meeting four lines in space.\label{F:four_lines}}
\end{figure}

The Schubert calculus~\cite{KL72} is a method to compute the number of solutions to
\Blue{{\it Schubert problems}}, which are a class of geometric problems
involving linear subspaces.
The prototypical Schubert problem is the classical problem of four lines:
How many lines in space meet four given lines?
To answer this, note that three lines $\ell_1,\ell_2,\ell_3$ lie on a unique doubly-ruled
hyperboloid, depicted in Figure~\ref{F:four_lines}.
These three lines lie in one ruling, while the second ruling consists of
the lines meeting the given three lines.
The fourth line $\ell_4$ meets the hyperboloid in two points.
Through each of these points there is a line in the second
ruling, and these are the two lines $m_1$ and $m_2$ meeting our four given lines.

The Galois group of this Schubert problem is the group of permutations
which are obtained by following the solutions over loops in the space of lines
$\ell_1,\dotsc,\ell_4$.
Rotating $\ell_4$ about the point $p$ gives a loop which interchanges the
two solution lines $m_1$ and $m_2$, showing that the Galois group
is $S_2$, the full symmetric group on two letters.

%
\subsection{Schubert problems in the Grassmannian}
%
A typical Schubert problem asks for the linear subspaces of a fixed dimension
(\Blue{{\it $k$-planes}}) in
$\C^n$ that have specified positions (incidence conditions) with respect to some
fixed, but otherwise general, linear subspaces.
Each incidence condition defines a set of $k$-planes, called a Schubert variety, and the
solutions to the Schubert problem are the points of intersection of the corresponding
Schubert varieties.
We describe this class of problems.

The \Blue{{\it Grassmannian $G(k,n)$}} is the set of $k$-planes in $\C^n$.
This is a complex manifold of dimension $k(n{-}k)$.
The problem of four lines involves the four-dimensional Grassmannian $G(2,4)$ as a line in
(projective) $3$-space corresponds to a 2-plane in $\C^4$.
The set of lines $m$ meeting a fixed line $\ell$ corresponds to
the set of 2-planes $M$ of $\C^4$ whose intersection with a fixed
2-plane $L$ is at least one-dimensional,
and this set of lines $m$ is a Schubert variety.
The problem of four lines asks for the points common to four such Schubert varieties,
one for each of the given lines $\ell_1$---$\ell_4$ in projective 3-space.

The specified positions of $k$-planes in Schubert problems are in reference to
flags in $\C^n$.
A \Blue{{\it flag} $F_\bullet$} is a sequence of linear subspaces
\[
  F_\bullet\ \colon\  F_1\subset F_2\subset\dotsb\subset
     F_{n-1}\subset F_n=\C^n\,,
\]
where $i=\dim F_i$.
The possible positions are encoded by partitions.
A \Blue{{\it partition} $\lambda$} is a weakly decreasing sequence of integers
\[
  \lambda\ \colon\ (n{-}k)\geq \lambda_1\geq \lambda_2\geq\dotsb\geq\lambda_k\geq 0\,.
\]
Give a partition $\lambda$ and a flag $F_\bullet$, the
\Blue{{\it Schubert (sub)variety $Y_\lambda F_\bullet$}} of $G(k,n)$ is
 \begin{equation}\label{Eq:Schub_Var}
   Y_\lambda F_\bullet\ :=\
   \{ E\in G(k,n)\mid \dim E\cap F_{n-k+i-\lambda_i}\geq i,\ i=1,\dots,k\}\,.
 \end{equation}
This has codimension $\Blue{|\lambda|}:=\lambda_1+\dotsb+\lambda_k$ in $G(k,n)$.
When $\lambda=\Box:=(1,0,\dotsc,0)$,
\[
   Y_\Box F_\bullet\ =\ \{E\in G(k,n)\mid \dim E\cap F_{n-k}\geq 1\}\,,
\]
as the other conditions are redundant.
We call $\Box$ a \Blue{{\it simple Schubert condition}} and
$Y_\Box F_\bullet$  a \Blue{{\it simple Schubert variety}}.
It depends only upon $F_{n-k}$, so we also write
$Y_\Box F_{n-k}$.
All four Schubert varieties in the problem of four lines are simple.

A \Blue{{\it Schubert problem}} is a list $(\lambda^1,\dotsc,\lambda^m)$ of partitions
with $|\lambda^1|+\dotsb+|\lambda^m|=k(n{-}k)$.
By Kleiman's Transversality Theorem~\cite{Kl74}, if
$F^1_\bullet, \dotsc, F^m_\bullet$ are general, then
the intersection
 \begin{equation}\label{Eq:Schubert_intersection}
   Y_{\lambda^1} F^1_\bullet \cap Y_{\lambda^2} F^2_\bullet
   \cap \dotsb\cap Y_{\lambda^m} F^m_\bullet
 \end{equation}
is transverse and consists of finitely many $k$-planes.
The number \Blue{$d(\lambda^1,\dotsc,\lambda^m)$} of $k$-planes may be computed using the
algorithms in the Schubert calculus (see~\cite{KL72} or~\cite{Fu97} or  
the Introduction to~\cite{HSS98}).
The problem of four lines is an instance of the
Schubert problem $(\Box,\Box,\Box,\Box)$ in $G(2,4)$ and our analysis shows that
$d(\Box,\Box,\Box,\Box)=2$.

We study
Schubert problems in which all except possibly two Schubert conditions are simple.
A \Blue{{\it simple Schubert problem}} on $G(k,n)$ is one of the form
\[
   \bigl(\lambda,\ \mu,\ \underbrace{\Box,\ \dotsc,\
      \Box}_{k(n-k)-|\lambda|-|\mu|}\bigr)\,, 
\]
where $\lambda,\mu$ are not necessarily equal to $\Box$.
We speak of the \Blue{{\it simple Schubert problem
$(\lambda,\mu)$ on $G(k,n)$}} (the $k(n-k)-|\lambda|-|\mu|$ simple conditions $\Box$ are
understood).

The primary reason for limiting our study to simple Schubert problems in this paper 
is that these are Schubert problems that are
complete intersections, and the off-the-shelf software that we use restricts us to
complete intersections.

%
\subsection{Galois groups of Schubert problems}\label{S:Galois}
%

According to Harris~\cite{Ha79}, Jordan~\cite{J1870} showed how intrinsic structures of
some enumerative problems could be understood in terms of Galois theory.
Harris took the opposite approach---computing Galois groups of enumerative problems to
expose the intrinsic structure of an enumerative problem.
He showed that many enumerative problems have Galois group equal to the full symmetric
group, demonstrating that these problems had no underlying structures.

Harris' theory relating Galois groups to monodromy groups begins with a map
$f\colon U\to V$ of degree $d$ between irreducible complex algebraic varieties $U$ and
$V$.
The function field $\C(U)$ of $U$ is a degree $d$ extension of the function field $\C(V)$
of $V$.
These fields may be embedded into the field $K$ of germs of meromorphic functions on a
disc around a regular value $v\in V$ of $f$.
If $L\subset K$ is the normalization in $K$ of the extension $\C(U)/\C(V)$, then the
Galois group $G=\Gal(L/\C(V))$ acts faithfully on the $d$ points $f^{-1}(v)$, and this
gives an embedding $G\hookrightarrow S_d$, where $S_d$ is the symmetric group of the fiber
$f^{-1}(v)$.

Replacing $U$ and $V$ by Zariski open subsets if necessary, we may assume that the
map $f\colon U\to V$ is a degree $d$ covering.
A loop in $V$ based at $v$ has $d$ lifts to $U$, one for each point in the fiber
$f^{-1}(v)$.
Associating a point in the fiber $f^{-1}(v)$ to the endpoint of the corresponding lift
gives a permutation in $S_d$.
This defines the usual permutation action of the fundamental group
of $V$ on the fiber $f^{-1}(v)$.
The \Blue{{\it monodromy group}} of the map $f\colon U\to V$ is the image of the
fundamental group of $V$ in $S_d$.

\begin{proposition}[Harris~\cite{Ha79}]\label{P:Galois}
  For a map $f\colon U\to V$ as above, the monodromy group equals the Galois
  group.
\end{proposition}

Given a Schubert problem $(\lambda^1,\dotsc,\lambda^m)$ on $G(k,n)$,
let $V$ be the space of $m$-tuples $(F_\bullet^1,\dotsc,F_\bullet^m)$ of flags,
a product of flag manifolds.
Define $U$ to be the incidence variety
 \begin{equation}\label{Eq:incidence}
    U\ :=\ \{(H,F_\bullet^1,\dotsc,F_\bullet^m)\in G(k,n)\times V\mid
     H\in Y_{\lambda^i}F_\bullet^i\ \mbox{\ for\ }i=1,\dotsc,m\}\ .
 \end{equation}
The fiber of $U$ over a point $H\in G(k,n)$ is a product of the Schubert subvarieties
\[
   \{F_\bullet^i \mid H\in Y_{\lambda^i}F_\bullet^i\}
   \qquad i=1,\dotsc,m
\]
of the flag manifold.
Each of these is irreducible, and so $U$ is irreducible.
Let $f\colon U\to V$ be the other projection.
Given $v=(F_\bullet^1,\dotsc,F_\bullet^m)\in V$, the fiber $f^{-1}(v)$
is the intersection~\eqref{Eq:Schubert_intersection}.
When $v$ is general, this has
$d=d(\lambda^1,\dotsc,\lambda^m)$ points, so that the map $f$ has degree $d$.
The \Blue{{\it Galois group}} of the Schubert problem is the Galois group
of the extension $\C(U)/\C(V)$.
By Proposition~\ref{P:Galois}, this is the 
monodromy group of the map $f\colon U\to V$.
\smallskip

The point of this paper is that these monodromy groups may be computed
using numerical homotopy continuation.
For this, we first compute the points in a single fiber $f^{-1}(v)$.
Then, given a loop $\varphi\colon[0,1]\to V$  based at $v$ ($\varphi(0)=\varphi(1)=v$), we
numerically follow the points in the fibers $f^{-1}(\varphi(t))$ as $t$ runs from $0$ to
$1$. 
This computes the lifts of $\varphi$ and thus the associated monodromy
permutation. 
Computing sufficiently many of these monodromy permutations will enable us to recover the
Galois group.
While this gives the idea behind our method, we postpone more details until
Section~\ref{subsec:loops}.

%
\section{Homotopy continuation of simple Schubert problems}\label{Sec:Homotopy}
%

Homotopy continuation is a numerical method to compute all solutions to a system of
polynomials given the solutions to a similar system.
We use it to find all solutions to a simple Schubert problem and to compute elements
of the monodromy group.
We first describe the method of numerical homotopy continuation, then discuss polynomial
formulations of Schubert problems, and finally explain the Pieri homotopy
algorithm~\cite{HSS98,HV00} to find all solutions to simple Schubert problems.

%
\subsection{Homotopy continuation of polynomial systems}\label{S:homotopy}
%

Suppose that we want to find all solutions to a 0-dimensional
\Blue{{\it target system}} of polynomial equations
 \begin{equation}\label{Eq:Poly_system}
   f_1(x_1,\dotsc,x_n)\ =\    f_2(x_1,\dotsc,x_n)\ =\ \dotsb\ =\
   f_N(x_1,\dotsc,x_n)\ =\ 0\,,
 \end{equation}
written as $F(x)=0$.
Numerical homotopy continuation finds these solutions if we have a
\Blue{{\it homotopy}}, which is a system \Blue{$H(x,t)$} of polynomials in $n{+}1$
variables such that
\begin{enumerate}
 \item The systems $H(x,1)=0$ and $F(x)=0$ both have the same solutions;

 \item We know all solutions to the \Blue{{\it start system}}
        $H(x,0)=0$;

 \item The components of the variety defined by $H(x,t)=0$ include curves whose projection
        to $\C$ (via the second coordinate $t$) is dominant; and

 \item The solutions to the system $H(x,t)=0$, where $t\in[0,1)$, occur at smooth points of
        curves from (3) in the variety $H(x,t)=0$.
\end{enumerate}

Given this, we restrict the variety $H(x,t)=0$ to $t\in[0,1]$ and obtain
finitely many real arcs in $\C^n\times[0,1]$ which connect (possibly singular) solutions
of the target system $H(x,1)=0$ to solutions of the start system $H(x,0)=0$.
We then numerically trace each arc from $t=0$ to $t=1$, obtaining all isolated
solutions to the target system.

\begin{figure}[htb]
\[
  \begin{picture}(121,130)(0,-25)
   \put(0,0){\includegraphics{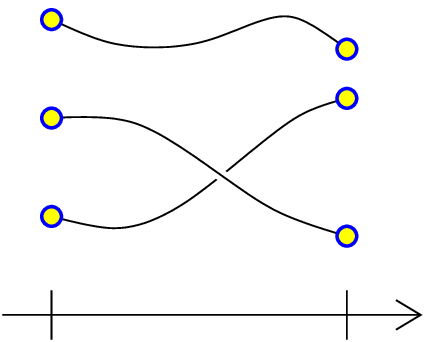}}
   \put(12,-11){$0$}   \put(97,-11){$1$}
   \put(108,14){$t$}
   \put(35,-25){optimal}
  \end{picture}
  \qquad
  \begin{picture}(121,130)(0,-25)
   \put(0,0){\includegraphics{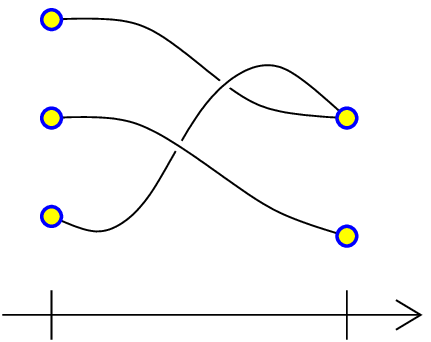}}
   \put(12,-11){$0$}   \put(97,-11){$1$}
   \put(28,-25){not optimal}
   \put(108,14){$t$}
  \end{picture}
  \qquad
  \begin{picture}(121,130)(0,-25)
   \put(0,0){\includegraphics{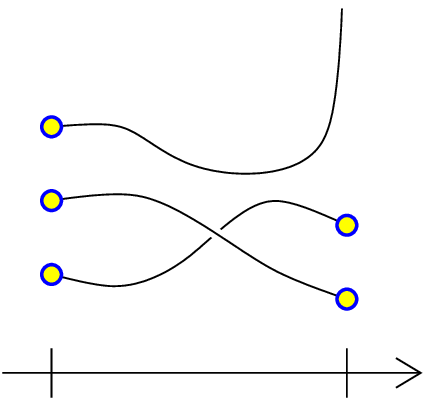}}
   \put(12,-11){$0$}   \put(97,-11){$1$}
   \put(28,-25){not optimal}
   \put(108,14){$t$}
  \end{picture}
\]
%
\caption{Optimal and non-optimal homotopies}
\label{F:homotopy_arcs}
\end{figure}

The homotopy is \Blue{{\it optimal}} if every solution at $t=0$ is connected to a unique
solution at $t=1$ along an arc.
This is illustrated in Figure~\ref{F:homotopy_arcs}.
For simple Schubert problems, the Pieri homotopy algorithm is optimal.

\begin{remark}
  \label{rem: actual homotopy}
  Homotopy continuation software often constructs a homotopy as follows.
  Let $F(x)$ be the \Blue{{\it target system}}~\eqref{Eq:Poly_system} and suppose we
  have solutions to a \Blue{{\it start system}} $G(x)$.
  Then for a number $\gamma\in\C$ with $|\gamma|=1$ define the
  \Blue{{\it linear homotopy}}
 \[
    H(x,t)\ :=\  \gamma t F(x) + (1-t) G(x)\,.
\]
  Then $H(x,t)$ satisfies the definition of a homotopy for all but
  finitely many $\gamma$.
  The software detects the probability 0 event that $H(x,t)$ does not satisfy the
  definition when it encounters a singularity, and
  then it recreates the homotopy with a different number $\gamma$.
\end{remark}

Path following algorithms use predictor-corrector methods,
which are conceptually simple for \Blue{{\it square systems}},
where the number of equations equals the number of variables.

Given a point $(x^{(0)},t^{(0)})$ on an arc such that $t^{(0)}\in[0,1)$, the $n\times n$
  matrix  
 \[
   H_x\ :=\ \left( \frac{\partial H_i}{\partial x_j}\right)_{i,j=1}^n
 \]
is regular at $(x^{(0)},t^{(0)})$, which follows from the definition of the homotopy. 
Let $H_t := (\partial H_1/\partial t, \dots, \partial H_n/\partial t)^T$.
Given $\Delta t$, we set
 \[
   \Blue{\Delta x}\ :=\  - \Delta t\, H_x(x^{(0)},t^{(0)})^{-1}\, H_t(x^{(0)},t^{(0)})\,.
 \]
For $t^{(1)} = t^{(0)} + \Delta t$, the point $(x',t^{(1)})=(x^{(0)}+\Delta x, t^{(1)})$ is
an approximation to the point $(x^{(1)},t^{(1)})$ on the same arc.
This constitutes a first order predictor step.
A corrector step uses the multivariate Newton method
for the system $H(x,t^{(1)})=0$, refining the approximate solution
$x'$ to a solution $x^{(1)}$.
In practice, the points $x^{(0)}$ and $x^{(1)}$ are numerical (approximate) solutions,
and both the prediction and correction steps require that $\det H_x \neq 0$
at every point where the computation of the Jacobian matrix $H_x$ is done.

When the system is not square,
additional strategies must be employed to enable the path following.
Fortunately, simple Schubert problems are exactly the class of Schubert problems for which
we have an optimal square homotopy (the Pieri homotopy).

\Blue{{\it Cheater homotopies}}~\cite{LSY89} are optimal homotopies
constructed from families of polynomial systems.
For example, given a Schubert problem $(\lambda^1,\dotsc,\lambda^m)$, let $V$ be the
space of all $m$-tuples $(F_\bullet^1,\dotsc,F_\bullet^m)$ of flags.
The total space of the Schubert problem
\[
    U\ :=\ \{(H,F_\bullet^1,\dotsc,F_\bullet^m)\in G(k,n)\times V\mid
     H\in Y_{\lambda^i}F_\bullet^i\ \mbox{ for } i=1,\dotsc,m\}
\]
is defined by equations (see Section~\ref{S:Equations}) depending upon the
point $(F_\bullet^1,\dotsc,F_\bullet^m)\in V$.
If $\varphi\colon\C\to V$ is an embedding of $\C$ into $V$ in which $\varphi(0)$ and
$\varphi(1)$ are general $m$-tuples of flags and we write
$\varphi(t)=(F_\bullet^1(t),\dotsc,F_\bullet^m(t))$, then
\[
  \varphi^* U\ =\
    \{(H,F_\bullet^1(t),\dotsc,F_\bullet^m(t)) \mid
     H\in Y_{\lambda^i}F_\bullet^i(t)\ \mbox{ for } i=1,\dotsc,m\}\ .
\]
This is defined by a system $H(x,t)=0$, which gives an optimal homotopy.
We use this particular cheater homotopy to compute permutations in monodromy groups,
called \Blue{{\it monodromy permutations}}.
We give more details in Section~\ref{subsec:loops}. 

There, we describe
the Pieri homotopy algorithm, 
which is a cheater homotopy where 
only one flag in $\varphi(t)$
actually moves and the others remain fixed.
The moving flag is in general position when $t=1$, but in a particular special position
when $t=0$, so that the Schubert problem
becomes a union of other Schubert problems (whose solutions were previously
computed and thus are known).

%
\subsection{Equations for Schubert problems}\label{S:Equations}
%
Polynomial homotopy continuation methods require that our geometric problems 
are modeled by a system of polynomial equations.
For efficiency, the number of variables should be minimized.
We describe equations for Schubert varieties and then model Schubert
problems by systems of equations which minimize the number of variables,
stated in Proposition~\ref{P:equations} below.

Represent a $k$-plane in $\C^n$ as the row space of
a $k$ by $n$ matrix $E$ with full rank and
a flag by an invertible $n$ by $n$ matrix $F_\bullet$ of constants,
where the $i$-dimensional subspace in the flag is the row space
of the first $i$ rows $F_i$ of the matrix.
The condition from~\eqref{Eq:Schub_Var} that $\dim E\cap F_{n-k+i-\lambda_i}\geq i$ is
 \begin{equation}\label{Eq:Schubert_equations}
   \mbox{rank} \left[\begin{matrix}E\\F_{n-k+i-\lambda_i}\end{matrix}\right]
   \ \leq\ n-\lambda_i\,,
 \end{equation}
which is given by the vanishing of the determinants of all
$n{+}1{-}\lambda_i$ by $n{+}1{-}\lambda_i$ submatrices of the
$n{+}i{-}\lambda_i$ by $n$ matrix in~\eqref{Eq:Schubert_equations}.
When $\lambda_i=0$ the condition~\eqref{Eq:Schubert_equations} is empty.

Write $\calE(E,F_\bullet,\lambda)$ for the system consisting of these
$\sum_i\binom{n+i-\lambda_i}{n+1-\lambda_i}\binom{n}{n+1-\lambda_i}$ equations.
The codimension $|\lambda|$ equals the number of equations only when $\lambda=\Box$.
In that case, $\calE(E,F_\bullet,\Box)$ consists of the single equation
 \begin{equation}\label{Eq:SiSchuCond}
    \det\ \left[\begin{matrix}E\\F_{n-k}\end{matrix}\right]\ =\ 0\,.
 \end{equation}

Since any two flags in general position are conjugate under a linear transformation,
we always assume that two flags
in~\eqref{Eq:Schubert_intersection} are fixed.
Let the flag $F_\bullet$ be defined by setting $F_i$ to be the span of
$e_n,e_{n-1},\dotsc,e_{n+1-i}$, where  $e_1,\dotsc,e_n$ form the standard basis
for $\C^n$.
The Schubert variety $Y_\lambda F_\bullet$ has an open subset isomorphic
to $\C^{k(n-k)-|\lambda|}$ consisting of $k$-planes that are the row space
of an echelon matrix of the form
\[
   \left[ \begin{matrix}
    0\dotsb 0 &\Blue{1}&\Magenta{*\dotsb *}&0&\Magenta{*\dotsb *}&0&\Magenta{*\dotsb *}\\
    0\dotsb 0 &0&0\dotsb 0  &\Blue{1}&\Magenta{*\dotsb *}&0&\Magenta{*\dotsb *}\vspace{-4pt}\\
               & \vdots &   &\vdots& \Blue{\ddots}  &\vdots& \Magenta{\vdots}\\
      0\dotsb 0 &0&0\dotsb 0&0&0\dotsb 0&\Blue{1}&\Magenta{*\dotsb *}
    \end{matrix}\right]
\]
where $1+\lambda_k, 2+\lambda_{k-1},\dotsc,k+\lambda_1$ are the columns with $1$s
and \Magenta{$*$} represents some  number.

To further reduce the number of variables, let the flag $F'_\bullet$ be defined by setting
$F'_i$ to be the span of $e_1,\dotsc,e_i$.
The \Blue{{\it skew Schubert variety}}~\cite{St77} (or Richardson variety),
\[
   \Blue{Y_{\lambda,\mu}}\ :=\
     Y_\lambda F_\bullet \cap Y_\mu F'_\bullet
\]
has an open subset parameterized by matrices of the form
\[
   E_{\lambda,\mu}\ :=\
   \left[\begin{array}{l}
     0\;\;\Blue{1}\;\;\Magenta{*\dotsb\dotsb\dotsb *}\hspace{4pt} 0\,\dotsb\, 0\vspace{-6pt}\\
     \,\vdots\hspace{1pt}\ddots\Blue{\ddots}\,\Magenta{*\dotsb\dotsb\dotsb *}
       \hspace{1pt}\ddots\hspace{1pt}\vdots\vspace{-1pt}\\
    0\,\dotsb\, 0\hspace{8pt}\Blue{1}\hspace{6pt}
     \Magenta{*\dotsb\dotsb\dotsb *}\hspace{3pt}0
   \end{array}\right]
\]
whose entries $a_{i,j}$ are
\[
   \left\{\begin{array}{rcl}
      1&\ &\mbox{if } j= i+\lambda_{k+1-i}\,,\\
      \Magenta{*}&&\mbox{if } i+\lambda_{k+1-i}<j\leq n-k+i-\mu_i\,,\\
      0&&\mbox{otherwise}\,.\end{array}\right.
\]
This parameterization is  one-to-one if the product of the rightmost
entries is non-zero,
\[
   0\ \neq\ \prod_{i=1}^k a_{i,\,n-k+i-\mu_i}\,.
\]

On the left below is $E_{\Box,\Box}$ in $G(2,4)$ and on the
right is $E_{210,110}$ in $G(3,7)$.
 \begin{equation}\label{Eq:loc_coords}
   \left[\begin{matrix}\Blue{1}&\Magenta{x}&0&0\\
                       0&0&\Blue{1}&\Magenta{y}
         \end{matrix}\right]
    \qquad\qquad
   \left[\begin{matrix}\Blue{1}&\Magenta{a}&\Magenta{b}&\Magenta{c}&0&0&0\\
                       0&0&\Blue{1}&\Magenta{d}&\Magenta{e}&0&0\\
                       0&0&0&0&\Blue{1}&\Magenta{f}&\Magenta{g}
          \end{matrix}\right]
 \end{equation}
%
%

Given a Schubert problem $\lambda,\mu,\nu^1,\dotsc,\nu^m$, we will always take two of the
flags to be these coordinate flags $F_\bullet$ and $F'_\bullet$, and
consider intersections of the form
\[
   Y_{\lambda,\mu}\cap Y_{\nu^1} F^1_\bullet \cap \dotsb\cap Y_{\nu^m} F^m_\bullet\,,
\]
where the flags $F^1_\bullet,\dotsc, F^m_\bullet$ are general.
By Kleiman's transversality theorem, all intersections
will lie in the subset of $Y_{\lambda,\mu}$ that is parameterized by matrices from
$E_{\lambda,\mu}$, and thus are solutions to the system of equations given by
\[
   \calE(E_{\lambda,\mu},F^1_\bullet,\nu^1),\ \dotsc,\
   \calE(E_{\lambda,\mu},F^m_\bullet,\nu^m)\,.
\]

This system is not necessarily square unless it is a
simple Schubert problem.
Since the homotopy continuation software we use is for square systems of polynomials, we
restrict ourselves to simple Schubert problems.
Write $G$ for the $n{-}k$ by $n$ matrix $F_{n-k}$.
Then $\calE(E_{\lambda,\mu},F_\bullet,\Box)$ is a single
equation~\eqref{Eq:SiSchuCond} that depends only on $G$.

\begin{proposition}\label{P:equations}
 A simple Schubert problem $(\lambda,\mu)$ on $G(k,n)$ is given by
 $\Blue{m}:=k(n-k)-|\lambda|-|\mu|$ matrices $G_1,\dotsc,G_m$ each of size
 $n{-}k$ by $n$, and the solutions are modeled by the system of equations
 \begin{equation}\label{Eq:SS_Eqs}
    \det\ \left[\begin{matrix}E_{\lambda,\mu}\\G_1\end{matrix}\right]\ =\
    \det\ \left[\begin{matrix}E_{\lambda,\mu}\\G_2\end{matrix}\right]
     \ =\ \dotsb \ =\
    \det\ \left[\begin{matrix}E_{\lambda,\mu}\\G_m\end{matrix}\right]\ =\ 0\,.
 \end{equation}
\end{proposition}

For example, the simple Schubert problem $\Box,\Box$ on $G(2,4)$ is modeled by
 \begin{equation}\label{Eq:G24}
   \det \left[\begin{matrix}\Blue{1}&\Magenta{x}&0&0\\ 0&0&\Blue{1}&\Magenta{y}\\
      g_{11}&g_{12}&g_{13}&g_{14}\\g_{21}&g_{22}&g_{23}&g_{24} \end{matrix}\right]
    \quad=\quad
   \det \left[\begin{matrix}\Blue{1}&\Magenta{x}&0&0\\ 0&0&\Blue{1}&\Magenta{y}\\
      g'_{11}&g'_{12}&g'_{13}&g'_{14}\\
      g'_{21}&g'_{22}&g'_{23}&g'_{24} \end{matrix}\right]
    \quad=\quad 0\,,
 \end{equation}
where $G_1=(g_{ij})$ and $G_2=(g'_{ij})$ are matrices of constants.

%
\subsection{Pieri homotopy algorithm}\label{S:Pieri}
%

We describe the simplified version of
the Pieri homotopy algorithm~\cite{HSS98,HV00}
that we use. 
The Pieri homotopy algorithm
 finds all solutions to
those Schubert problems where all except possibly two  partitions
consist of a single 
part, $(a,0,\dotsc,0)$. 
It is based on subtle geometric degenerations constructed in~\cite{So97}.
Both the algorithm and the degenerations enjoy a dramatic simplification for simple
Schubert problems.
The degenerations 
for these simple Schubert problems
were introduced by Schubert~\cite{Sch1886c,Sch1886b}.

\begin{example}\label{Ex:Pieri_homotopy}
 Consider the simple Schubert problem $(\lambda,\mu)=(210,110)$ in $G(3,7)$.
 Begin with local coordinates~\eqref{Eq:loc_coords} for $E_{210,110}$
 \[
   E\ :=\ E_{210,110}\ =\
   \left[\begin{matrix}\Blue{1}&\Magenta{a}&\Magenta{b}&\Magenta{c}&0&0&0\\
                       0&0&\Blue{1}&\Magenta{d}&\Magenta{e}&0&0\\
                       0&0&0&0&\Blue{1}&\Magenta{f}&\Magenta{g}
          \end{matrix}\right]
\]
 There are four columns, $1$, $2$, $3$, and $6$ not of the form
 $n{+}i{-}k{-}\mu_i$.
 Let $G$ be a general $4$-plane represented by a matrix in which these columns
 form a identity matrix
 \begin{equation}\label{Eq:gen_G}
  G\ :=\ \left[
    \begin{matrix}\Blue{1}&0&0&\Magenta{*}&\Magenta{*}&0&\Magenta{*}\\
                  0&\Blue{1}&0&\Magenta{*}&\Magenta{*}&0&\Magenta{*}\\
                  0&0&\Blue{1}&\Magenta{*}&\Magenta{*}&0&\Magenta{*}\\
                  0&0&0&\Magenta{*}&\Magenta{*}&\Blue{1}&\Magenta{*}\end{matrix}\right]\ .
 \end{equation}
 Let $G(t)$ be this matrix with each entry $\Magenta{*}$ scaled by $t$.
 Then
 \begin{equation}\label{Eq:Pieri_H_equation}
     \det\ \left[\begin{matrix} E\\G(t)\end{matrix}\right]\ =\
     ceg + t(\Magenta{*}c+\dotsb+\Magenta{*}cef) +
      t^2(\Magenta{*}+\dotsb+\Magenta{*}bef)    +t^3(\Magenta{*}f+\Magenta{*}af)\,,
 \end{equation}
%
where each \Magenta{$*$} again represents a fixed number.

When $t=0$, the expression~\eqref{Eq:Pieri_H_equation} becomes $ceg$.
Let us investigate the consequences of $ceg=0$.
If we set $c=0$ in $E_{210,110}$, we get $E_{210,210}$ and if we set $g=0$, we get
$E_{210,111}$,
 \[
   E_{210,210}\ =\
   \left[\begin{matrix}\Blue{1}&\Magenta{a}&\Magenta{b}&\Brown{0}&0&0&0\\
                       0&0&\Blue{1}&\Magenta{d}&\Magenta{e}&0&0\\
                       0&0&0&0&\Blue{1}&\Magenta{f}&\Magenta{g}
          \end{matrix}\right]
  \qquad
   E_{210,111}\ =\
   \left[\begin{matrix}\Blue{1}&\Magenta{a}&\Magenta{b}&\Magenta{c}&0&0&0\\
                       0&0&\Blue{1}&\Magenta{d}&\Magenta{e}&0&0\\
                       0&0&0&0&\Blue{1}&\Magenta{f}&\Brown{0}
          \end{matrix}\right]\ .
\]
(If $e=0$, then the row operation $R_1\leftarrow R_1-cR_2$ gives a matrix with $c=0$,
 which lies in $E_{210,210}$.)
This computation in local coordinates shows that
\[
  Y_{210,110}\cap Y_\Box G(0)\ =\
  Y_{210,210}\cup Y_{210,110}\,.
\]

Now suppose that we have general $4$ by $7$ matrices $G_1,\dotsc,G_6,G_7$,
and we wish to solve the instance of the simple Schubert problem
$(210,110)$:
 \begin{equation}\label{Eq:G37}
    \det\ \left[\begin{matrix}E\\G_1\end{matrix}\right]\ =\
    \dotsb \ =\
    \det\ \left[\begin{matrix}E\\G_6\end{matrix}\right]\ =\ 
    \det\ \left[\begin{matrix}E\\G_7\end{matrix}\right]\ =\ 0\,,
 \end{equation}
where $G_7$ is the matrix~\eqref{Eq:gen_G}.
Replacing $G_7$ by  $G(t)$ gives a homotopy of 7 equations in the coordinates
$a,\dotsc,g$ where 6 equations are fixed~\eqref{Eq:G37} and one depends on
$t$~\eqref{Eq:Pieri_H_equation}.
When $t=0$, the latter becomes $ceg=0$ and the system splits into subsystems
on $E_{210,210}$ and $E_{210,111}$ involving the matrices $G_1,\dotsc,G_6$.
Numerical continuation along this homotopy uses solutions to these smaller problems to
obtain solutions to the system~\eqref{Eq:G37}.
\end{example}

Write $\lambda\lessdot\nu$ if the components of the vector $\nu-\lambda$ are either 0 or
1, with exactly one~1.
For example, $110\lessdot 210$ and $110\lessdot 111$.
Given partitions $\lambda,\mu$, define the $(n{-}k)$-plane 
 \[
   G_\mu \ :=\ \{ e_i\mid i\not\in\{n{-}k{+}j{-}\mu_j\,,\ \mbox{for}\ j=1,\dotsc,k\}\}\,.
 \]
Then if the non-zero entries of the matrix $E_{\lambda,\mu}$ are $a_{i,j}$, we have
 \begin{equation}\label{Eq:Degenerte_eqs}
   \det\ \left[\begin{matrix}E_{\lambda,\mu}\\G_\mu\end{matrix}\right]\ =\
      \prod_{i=1}^k a_{i,n-k+i-\mu_i}\ ,
 \end{equation}
%
which is the product of the rightmost non-zero entries in
the rows of $E_{\lambda,\mu}$.
%
This determinant defines the charts $E_{\lambda,\nu}$ for $\mu\lessdot \nu$.
A \Blue{{\it child problem}} for the simple Schubert problem $(\lambda,\mu)$
is one of the form 
$(\lambda,\nu)$ with $\mu\lessdot\nu$.

The Pieri homotopy algorithm finds all solutions to a simple Schubert
problem  $(\lambda,\mu)$ with fixed (but general) $(n{-}k)$-planes $G_1,\dotsc,G_m$
where $m+|\lambda|+|\mu|=k(n{-}k)$.
We assume that we are given all solutions to all 
child problems $(\lambda,\nu)$ with $\mu\lessdot\nu$
and the $(n{-}k)$-planes are  $G_1,\dotsc,G_{m-1}$.
If we let $G_m(t)$ be a 1-parameter family of  $(n{-}k)$-planes with $G_m(1)=G_m$ and
$G_m(0)=G_\mu$, then we obtain a homotopy
 \[
    \det\ \left[\begin{matrix}E_{\lambda,\mu}\\G_1\end{matrix}\right]
     \ =\ \dotsb \ =\
    \det\ \left[\begin{matrix}E_{\lambda,\mu}\\G_{m-1}\end{matrix}\right]\ =\
    \det\ \left[\begin{matrix}E_{\lambda,\mu}\\G_m(t)\end{matrix}\right]\ =\ 0\,.
 \]
When $t=0$ this is the disjunction of child problems and when $t=1$, it is
the problem we wish to solve.

This method recursively finds solutions to parent problems given solutions to their
child problems.
The depth of this recursion equals the dimension of the skew Schubert variety
corresponding to the simple Schubert problem we wish to solve, 
i.e., the number of variables in the corresponding equations.
The base case of this recursion
is when $|\lambda|+|\mu|=k(n{-}k)$, for then
$E_{\lambda,\mu}$ is empty unless $\lambda_i+\mu_{k+1-i}=n{-}k$ for $i=1,\dotsc,k$, and in
that case, $E_{\lambda,\mu}$ gives the $k$-plane spanned by
$\{e_{i+\lambda_i}\mid i=1,\dotsc,k\}$.

These homotopies are optimal.
This is because they only follow solutions to the given Schubert problem and because the
number $d(\lambda,\mu)$ of solutions to a simple Schubert problem
satisfies the same recursion as the number of paths that are followed.
Namely,  if $|\lambda|+|\mu|=k(n{-}k)$ then  $d(\lambda,\mu)=0$ unless
$\lambda_i+\mu_{k+1-i}=n{-}k$ for $i-1,\dotsc,k$, and then it equals 1.
If $|\lambda|+|\mu|< k(n{-}k)$, then
 \[
   d(\lambda,\mu) 
   \quad = \quad \sum_{\mu\lessdot\nu} d(\lambda,\nu)\ .
 \]
%

\begin{remark}
 We do not quite use the homotopy we just described, as the equations involving $t$ will
 in general have degree in $t$ at least the minimum of $k$ and $n{-}k$.
 Instead, we use the convex combination of the equations
 \begin{equation}\label{Eq:pencil}
   \gamma
   t \det \left[\begin{matrix}E_{\lambda,\mu}\\G_m\end{matrix}\right]
   \ +\ (1-t)\det\left[\begin{matrix}E_{\lambda,\mu}\\G_\mu\end{matrix}\right]
   \ =\ 0\,.
 \end{equation}
  Here, $\gamma \in \C$ has norm 1, $|\gamma|=1$.
 This has degree 1 in the homotopy parameter $t$.

 Doing this, the intermediate solutions will not necessarily be solutions to the Schubert
 problem, so we need to argue that this homotopy will remain optimal.
 For this, we appeal a little to the geometry of the Grassmannian.
 The equations we use define hyperplane sections of the Grassmannian in its Pl\"ucker
 embedding, and the number of solutions $d(\lambda,\mu)$ to the Schubert problem turns out 
 to be the degree of the variety $Y_{\lambda,\mu}$.
 Thus, replacing a family of hyperplanes defined by Schubert conditions (as in the Pieri
 homotopy) by an arbitrary pencil of hyperplanes~\eqref{Eq:pencil} will still give an
 optimal homotopy between solutions to the child problems and solutions to the parent
 problem. 
\end{remark}

%
\section{Description of software}
%

We divide the computation of the Galois group of
a Schubert problem into three tasks.
 \begin{enumerate}
  \item
    Compute the solutions of a general instance of the problem, called
    a \Blue{{\it master set}}.
  \item Use cheater homotopies to compute \Blue{{\it monodromy permutations}} of the
  master set.
\item Determine the \Blue{{\it group}} generated by the monodromy permutations.
\end{enumerate}

The first task may be accomplished by the brute-force application
of a polynomial system solver.
This is, however, inefficient.
In Table~\ref{table:Schub_Kouch}, we compare the number of solutions to simple Schubert
problems with $\lambda=\mu=\Box$ to the number of homotopy paths followed in
polyhedral homotopies (as in the black-box solver of PHCpack).
This is the volume of the associated Newton polytope and was computed with
Polymake~\cite{Polymake}. 
\begin{table}[htb]
  \begin{tabular}{|c||c|c|c|c|c|c||c|c|c||c|c||c|}\hline
  $k,n$ &2,6&2,7&2,8&2,9&2,10&2,11&3,6&3,7&\Blue{3,8}&4,6&4,7&5,7\\\hline
  solutions&14&42&132&429&1430&4862&42&462&\Blue{6006}&14&462&42\\\hline
  paths&18&67&248&919&3426&12843&130&3004&\Blue{74645}&42&7156&364\\\hline
 \end{tabular}\vspace{4pt}
 \caption{Inefficiency of polyhedral homotopy for simple Schubert problems}
 \label{table:Schub_Kouch}\vspace{-10pt}
\end{table}
This shows that the Pieri homotopy algorithm is an efficient alternative.

Subsection \ref{subsec:loops} discusses the second task.

The theory for the third task is beyond the scope of this paper.
Based on preliminary computations, we conjectured that the Galois group is always the
full symmetric group, and therefore we only check this.
The fastest routine we have found to accomplish this is
the {\tt isNaturalSymmetricGroup} function of GAP~\cite{GAP4}.

In Section~\ref{Sec:software}, we discuss implementations of our algorithm.
Implementations and documentation of our computations are available at our
website~\cite{LeSo-WWW}.

%
\subsection{Computing monodromy permutations}\label{subsec:loops}
%

Suppose that we have a master set of solutions to the
simple Schubert problem $(\lambda,\mu)$ on $G(k,n)$, which is
modeled by the system of equations   
 \begin{equation}\label{Eq:static_system}
  \det\ \left[\begin{matrix}E_{\lambda,\mu}\\G_1\end{matrix}\right]\ =\
  \det\ \left[\begin{matrix}E_{\lambda,\mu}\\G_2\end{matrix}\right]
   \ =\ \dotsb \ =\
  \det\ \left[\begin{matrix}E_{\lambda,\mu}\\G_m\end{matrix}\right]\ =\ 0\,,
 \end{equation}
where $G_1,\dotsc,G_m$ are fixed $(n{-}k)$-planes.
We follow these solutions along loops in the
space of $m$-tuples of $(n{-}k)$-planes to compute monodromy permutations.

The off-the-shelf homotopy continuation software we use requires that the equations are
linear in the homotopy parameter $t$, and so we follow piece-wise linear loops.
For these, we fix all of the $G_i$ except $G_m$, and exploit the
linearity of the determinant in the row vectors
$g_1,\dotsc,g_{n-k}$ of $G_m$.
If  $g'_i$ is a vector not in $G_m$ and we replace the row $g_i$ by the convex combination
$(1-t)g_i+tg'_i$, obtaining the pencil of planes $G_m(t)$, then
 \begin{eqnarray*}
   \det\ \left[\begin{matrix}E_{\lambda,\mu}\\G_m(t)\end{matrix}\right]
    &=&
    (1-t) \det\ \left[\begin{matrix}E_{\lambda,\mu}\\G_m\end{matrix}\right]
    \ +\ t \det\ \left[\begin{matrix}E_{\lambda,\mu}\\G_m(1)\end{matrix}\right]\\
    &=& (1-t)F(g_1,\dotsc,g_m)\ +\ tF(g_1,\dotsc,g'_i,\dotsc,g_m)\,.
 \end{eqnarray*}
Replacing the equation in~\eqref{Eq:static_system} involving $G_m$ with this
equation gives a linear homotopy between
the system~\eqref{Eq:static_system} and one with $G_m(1)$ in place of $G_m$.

Given a different $(n{-}k)$-plane $G'_m$ spanned by
$g'_1,\dotsc,g'_{n-k}$, we use these vectors to generate
loops along which we can compute monodromy permutations.
Suppose for illustration that $n{-}k=3$ and the vectors are $[a,b,c]$ for $G_m$ and
$[a',b',c']$ for $G'_m$.
The different pencils that we may create correspond to the edges of a cube.
\[
  \begin{picture}(166,118)(-83,-4)
                  \put(-20,102){$[a',b',c']$}
   \put( -20, 97){\line(-2,-1){30}}   \put( 20, 97){\line( 2,-1){30}}
                  \put(0,97){\line(0,-1){15}}
   \put( -82, 69){$[a',b',c]$}        \put( 43,  69){$[a,b',c']$}
                \put( -20,  69){$[a',b,c']$}
   \put( -60, 43){\line( 0,1){20}} 
   \put(-50, 43){\line( 2,1){18}}   \put(-28, 54){\line( 2,1){18}}
   \put(  60, 43){\line( 0,1){20}}  
   \put( 50, 43){\line(-2,1){18}}    \put( 28, 54){\line(-2,1){18}}
   \put( -10, 43){\line(-2,1){40}}  \put( 10, 43){\line( 2,1){40}}

   \put( -80, 31){$[a',b,c]$}        \put( 43,  31){$[a,b,c']$}
                \put( -17.2,  31){$[a,b',c]$}
   \put( -20, 10){\line(-2,1){30}}   \put( 20, 10){\line(2,1){30}}
                  \put(0,10){\line(0,1){15}}
                 \put( -17.2,  -2){$[a,b,c]$}
  \end{picture}
\]

Our software offers 3 strategies to generate loops.

\begin{itemize}

\item {\bf long loop} goes from a vertex of the cube to its
opposite and back,
\begin{eqnarray*}
[a,b,c] &\to\  [a',b,c]\ \to\ [a',b',c]\  \to & [a',b',c']\ \to\\
        &\to\  [a,b',c']\ \to\ [a,b,c']\ \to & [a,b,c]\,.
\end{eqnarray*}

\item {\bf short loop} uses only a square, 
\[
   [a,b,c]\ \to\ [a',b,c]\ \to\ [a',b',c]\ \to\
   [a,b',c]\ \to\ [a,b,c]\,.
\]

\item {\bf half loop} makes use of just one edge,
 \[
    [a,b,c]\ \stackrel{1}{\to} \ [a',b,c]\ \stackrel{\gamma}{\to}\ [a,b,c]\,,
 \]
where the second homotopy is modified via a random number $\gamma\in\C$,
\[
   H_\gamma\ :=\  (1-t)F(a',b,c)\ +\ \gamma t F(a,b,c)\,.
\]
\end{itemize}

While we do not offer a proof that these loops will suffice to find all
non-trivial permutations, we remark that they do suffice in the examples we considered.

\begin{figure}[htb]

\begin{picture}(400,230)(0,-10)
 \put(0,33.1){\begin{picture}(160,170)
   \put(0,0){
   \includegraphics{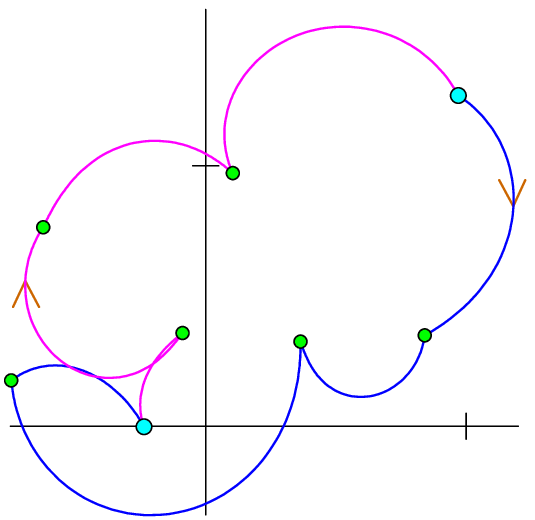}
   }
   \put(129,10){$1$} \put(45,97){$i$}
   \put(33,15){$m_1$}   \put(137,125){$m_2$}
   \put(0,-13){First coordinate in position $1,2$}
 \end{picture}}

 \put(200,0){\begin{picture}(160,220)
   \put(0,0){
   \includegraphics{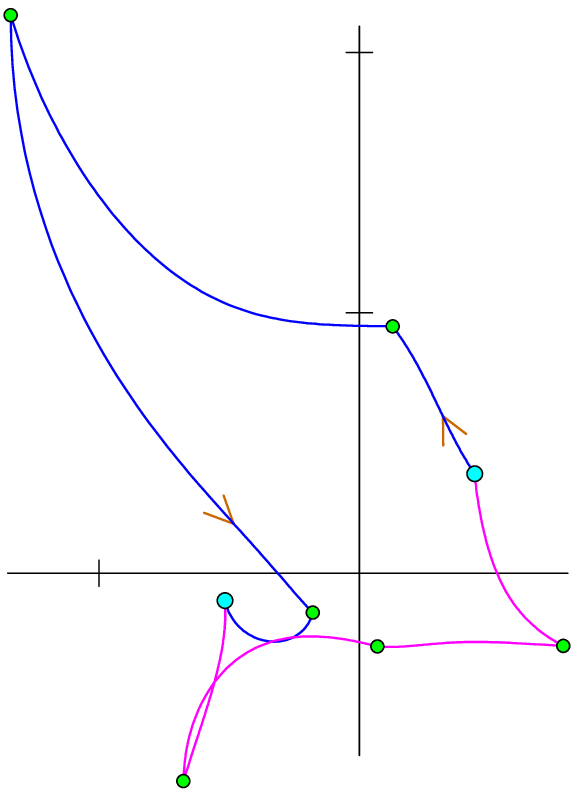}
   }
   \put(14,40){$-1$} \put(90,134){$i$}\put(85,207){$2i$}
   \put(47,48){$m_1$}   \put(143,85){$m_2$}
   \put(0,-13){Second coordinate in position $2,4$}
 \end{picture}}

\end{picture}

\caption{Paths tracked in Example~\ref{Ex:short_loops}.}\label{F:Monodromy_permutation}
\end{figure}

\begin{example}\label{Ex:short_loops}
 Suppose that we have the simple Schubert problem
 $(\Box,\Box)$ on $G(2,4)$ as given by~\eqref{Eq:G24} with
 \begin{eqnarray*}
  G_1& =& \left[\begin{matrix}
      -55 -  8 i &  17 + 15 i &  40 + 99 i & -17 - 38 i \\
      -67 + 25 i & -82 - 55 i & -99 - 80 i & -21 - 85 i\end{matrix}\right]
    \\
  G_2&=&\left[\begin{matrix}
     66 + 53 i & -73 - 14 i & 85 +  5 i &  67 + 16 i \\
    -53 - 85 i &  36 - 25 i &  2 + 81 i & -58 + 35 i \end{matrix}\right]\ .
 \end{eqnarray*}
 Its solutions $m_1$ and $m_2$ are
 \begin{eqnarray*}
  m_1&:=&\left[\begin{matrix}
    1&-0.23714-.0028980i &\ 0\ &0\\ 0&0&1&-.51680-.10520i \end{matrix}\right]
   \\
   m_2&:=&\left[\begin{matrix}
    1&.97009+1.2705i &\ 0\ &0\\ 0&0&1& .44336+.38248i\end{matrix}\right]
 \end{eqnarray*}
 The short loop strategy with 
\[
  G'\ =\ \left[\begin{array}{rrrr}
     33 - 84 i &  21 -    i & 59 + 94 i & -94 + 89 i\\
    -15 - 19 i &       29 i & 79 + 51 i &  89 + 3 i \end{array}\right]
\]
 creates a non-trivial monodromy permutation.
 The paths followed during the homotopies are drawn in
 Figure~\ref{F:Monodromy_permutation}, where the large circles are
 the values taken at the endpoints of the different homotopies.
\end{example}

%
\subsection{Implementation}\label{Sec:software}
%
We have two Maple implementations of our algorithms using the package
PHCmaple~\cite{PHCmaple-ICMS-06} to interface with
PHCpack~\cite{V99}, which performs the numerical polynomial homotopy continuation.
PHCmaple produced the graphic of Figure~\ref{F:Monodromy_permutation}.
The second implementation may alternatively call Bertini~\cite{Bertini}.

Our prototype implementation was carried out entirely in Maple to
take advantage of Maple packages to generate the equations and to manage the
monodromy group, while using the black-box solver in PHCpack to compute the master sets of 
solutions.
The largest problem this implementation could treat was the simple Schubert problem
$(210,200)$ on $G(3,7)$---it showed that the Galois group is the full symmetric group
$S_{91}$. 
Previously, the largest Schubert problem
whose Galois group that was 
proven
to be the full symmetric group was the simple Schubert problem $(20,10)$ on
$G(2,6)$ with 9 solutions~\cite{BV05}.

Our second implementation also uses Maple and either PHCpack or Bertini. 
However, it relies on the Pieri homotopy algorithm to compute master sets of solutions and GAP
to manage the monodromy groups, removing the two main computational bottlenecks of the 
prototype.
The largest problems this implementation has treated are 
$(2100,\Box)$ in $G(4,8)$ with 8580 solutions, 
$(210,210)$ in $G(3,8)$ with 10329 solutions and $(210,200)$ on $G(3,9)$ with 
\Blue{17589} solutions.
It showed that all of these have Galois group equal to the full symmetric group.
The computation for the problem $(\Box,\Box)$ on $G(3,8)$ with 6006 solutions is the basis
of the Numerical Theorem. 

We computed Galois groups of the simple Schubert problems
$(\Box,\Box)$ on all small Grassmannians, using the short loops strategy.
They were run on several different computers, including an 
AMD Athlon 64 Dual Core Processor 4600+ with CPU clock speed of
2400 MHz and 1 GB of memory whose timings (using PHCpack) are reported in
Table~\ref{table:timings}. 
\begin{table}[htb]
 \begin{tabular}{|c||c|c|c|c|c|c|c|}\hline
  $k,n$ &2,4&2,5&2,6&2,7&2,8&2,9&2,10\\\hline
  solutions&2&5&14&42&132&429&1430\\\hline
  time&12s&27s&19s&51s&4.2m&20.5m&2.6h\\\hline
  permutations&4&6&5&6&7&4&7\\\hline
 \end{tabular}\vspace{8pt}
 \begin{tabular}{|c||c|c|c|c|c||c|c|c|}\hline
  $k,n$ &3,5&3,6&3,7&3,8&\Blue{3,9}&4,6&4,7&4,8\\\hline
  solutions&5&42&462&6006&\Blue{17589}&14&462&8580\\\hline
  time&12s&35s&17.9m&18.6h&\Blue{78.2h}&15s&23.5m&44.5h\\\hline
  permutations&4&4&5&6&\Blue{7}&5&5&7\\\hline
 \end{tabular}
\begin{center} s := seconds, m := minutes, h := hours\end{center}
 \caption{Timings of Galois group computation}
 \label{table:timings}\vspace{-10pt}
\end{table}
These reported times are not CPU times, but actual elapsed (wall clock) time, and so may exceed
CPU time by 10 to 20 \%.
We also record the number of permutations we needed to compute.
The entry in $G(3,9)$ is the Schubert problem $(210,200)$ and the entry in $G(4,8)$ is the
Schubert problem $(2100,\Box)$. 

We ran some of these computations in Bertini.
It was unable to compute examples in more than 10 variables, and was markedly slower for
the largest computations it completed, on $G(2,8)$, $G(3,7)$, and $G(4,7)$.
On the other hand, Bertini provided an independent verification that the Galois groups were
indeed the full symmetric groups. 
We remark that Bertini is new software and its efficiency will likely improve.
%

%
\section{Conclusions and Future work}
%

This paper demonstrates the feasibility of homotopy continuation as tool to study the Galois 
groups of enumerative problems.
It also implicitly provides several challenges to the numerical homotopy community.
Perhaps the most serious is the current lack of certifiability of computations in numerical homotopy software.
While numerical methods will increasingly outperform symbolic algorithms in algebraic geometry, they are 
currently inferior in that their results do not come with certification.
Certificates for numerical computations do exist in theory, for example in Shub and Smale's~\cite{SS93} alpha-theory,
and there is a need for their implementation. In fact, even more reliable numerical techniques, 
such as the interval step control proposed in \cite{KX94}, are sidestepped by homotopy continuation 
software developers mostly due to the perceived complexity of implementation 
and the expected slower performance in comparison with heuristic methods.
Robust, off-the-shelf software to handle polynomial systems that are not complete intersections 
is also needed to deal with 
ideals in algebraic geometry, which are typically not complete intersections.

Two further theoretical problems are not addressed in this paper.
While the sampling of the fundamental group 
of the base space provided by the loop-generating heuristics of
subsection \ref{subsec:loops} 
is sufficient to 
generate the Galois group in the examples we considered, 
a better understanding of the topology of a complement of an algebraic variety in a Grassmannian 
is needed to prove that this sampling is always sufficient.
Second, 
we do not know how to certify 
that the computed set of permutations generates the whole Galois
group, if 
it is not the full symmetric group.

The tools that we use could be applied more systematically to other problems in enumerative geometry.
To this end,
we plan a comprehensive project 
exploring the limit of computability of Galois groups of
Schubert problems along several fronts.
This will involve the software and algorithms described here---perhaps also incorporating HOM4PS~\cite{HOM4PS}.
We will also write parallel software implementing algorithms to compute Galois groups of Schubert 
problems that are not compete intersections, as well as pushing the limits of the symbolic methods of Billey
and Vakil~\cite{BV05} and of Vakil's combinatorial algorithm~\cite{Va06b}.
This will be a large and distributed computation as in~\cite{RSSS}
which should give a catalog of the Galois groups of several tens of thousands of Schubert problems.

\providecommand{\bysame}{\leavevmode\hbox to3em{\hrulefill}\thinspace}
\providecommand{\MR}{\relax\ifhmode\unskip\space\fi MR }
\providecommand{\MRhref}[2]{%
  \href{http://www.ams.org/mathscinet-getitem?mr=#1}{#2}
}
\providecommand{\href}[2]{#2}


\end{document}